\documentclass[a4paper,11pt,twoside,leqno]{article}

\usepackage{amsmath,amsthm,amsfonts,latexsym,amscd,amssymb}
\usepackage[a4paper, total={6in,8in}]{geometry}
\numberwithin{equation}{section}
\input xypic
\def\qed{{\hbadness=10000\hfill\ \vbox{\hrule height.09ex
			\hbox{\vrule width.09ex height1.55ex depth.2ex \kern1.8ex
				\vrule width.09ex height1.55ex depth.2ex}\hrule height.09ex}\break
		\bigskip}}
\newtheorem{theorem}{Theorem}[section]

\newtheorem{corollary}{Corollary}[section]
\newtheorem{proposition}{Proposition}[section]

\theoremstyle{definition}
\newtheorem{definition}[theorem]{Definition}

\theoremstyle{remark}

\newcommand{\n}{\noindent}

\begin{document}
	
	\linespread{1}\title {\textbf{Relatively normal-slant helices in Minkowski $3$-space }}
	
\author{Akhilesh Yadav, Ajay Kumar Yadav\thanks{Corresponding author}\\{\n E-mail: ajaykumar74088@gmail.com\thanks{Corresponding author E-mail}}}
	\date{}
 
\maketitle 

	\noindent\textbf{Abstract:} In this paper, we study relatively normal-slant helices lying on timelike as well as spacelike surfaces in Minkowski $3$-space $ \mathbb{E}_1^3$. The axes of spacelike and timelike relatively normal-slant helices are obtained via their Darboux frames. We also establish characterization theorems for spacelike and timelike relatively normal-slant helices in Minkowski $3$-space $\mathbb{E}_1^3$. Finally, the relationship between relatively normal-slant helices and slant helices is found on timelike as well as spacelike surfaces.

\n\textbf {Mathematics Subject Classification (2020):}  53A04, 53A05, 53A35

\n\textbf{Key words:} Relatively normal-slant helices, Frenet-frame, Darboux-frame, asymptotic curve, geodesic curve, slant helix.

\section{Introduction}
 In recent years helical curves are studied widely in Euclidean and non-Euclidean spaces. Generalized helices, slant helices, isophote curves and relatively normal slant helices are such curves. When we study about curves, we need a frame of reference along the curve, as a result many moving frames along a curve have been established, among which Frenet frame and Darboux frame are prominent. We get a type of helical curve when one of the vector fields in the chosen frame makes a constant angle with a fixed vector. Helical curves are most fascinating curves that attracted attention in a wide range of disciplines such as Mathematics, Physics, Biology, Architecture and Computer Science.

\n In 3-dimensional Euclidean space $E^3$, a curve $\gamma$ with moving Frenet frame $\{ T, n, b\}$ is a helix if its tangent vector field $T$ makes a constant angle with a fixed direction $d$, called the axis of the helix. Helices are characterized by the fact that the ratio $\tau/\kappa$ is constant along the curve, where   $\tau$ and  $\kappa$  stands for the torsion and curvature of $\gamma$, respectively [3]. The notion of helix in Minkowski 3-space $E_1^3$ is developed similarly. In $[7]$, Izumiya and Takeuchi defined a slant helix in $E^3$, by the property that the principal normal $n$ to the curve makes a constant angle with a fixed direction, called its axis. In $[1]$, Ali and Lopez  looked into slant helices in Minkowski $3$-space $E_1^3$ and they obtained some characterizations of slant helix and its axis.

\n Using Darboux frame $ \{ T, B = N \times T, N\}$ of a curve, lying on a surface with unit normal $N$, In 2015, Dogan and Yayli investigated isophote curve, in the Euclidean space $[4]$. They found the axis of an isophote curve via its Darboux frame and then gave some characterizations about the isophote curve and its axis in Euclidean 3-space. Isophote curve is defined as a locus of points of a surface at which the normal to the surface makes a constant angle with a fixed direction, called its axis. Isophote curves are one of the characteristics curves on a surface such as geodesic, asymptotic or line of curvature. In [5], Dogan defined Isophote curves on timelike surfaces in Minkowski 3-space and found the axes of spacelike and timelike isophote curves via their Darboux frames. In [6], Dogan and Yayli defined Isophote curves on spacelike surfaces in Lorentz-Minkowski space and they found its axis as timelike and spacelike vectors via the Darboux frame. They also gave some relations between isophote curves and special curves on surfaces such as geodesic curves, asymptotic curves or lines of curvature.

\n Recently, Macit and Duldul $[9]$ defined the relatively normal-slant helix on a surface in the Euclidean space $E^3$ via  Darboux frame $\{T, B = N \times T, N\}$, along the curve whose vector field $B$ makes a constant angle with a fixed direction, called its axis. They gave some characterizations for such curves and obtained relations between some special curves (general helices, integral curves, etc.) and relatively
normal-slant helices. In this paper we study the relatively normal-slant helix in Minkowski 3-space $E_1^3$, where we obtain the axis of the relatively normal-slant helix and get some characteristics of the curve. Dogan $[5]$ while investigating isophote curves noticed that the curve which is both a geodesic and a slant helix is an isophote curve, interestingly we come to know that the curve which is both an asymptotic curve and a slant helix is a relatively normal-slant helix. The paper is arrange as follws: In section 2, we discuss some basic theory of unit speed parametrized curve on a smooth surface in Minkowski $3$-space $E_1^3$. In sections (3, 4, 5), we find the axis of a relatively normal-slant helix on a spacelike and timelike surface. In section 6, we establish characterization theorems for spacelike and timelike relatively normal-slant helices in Minkowski $3$-space $\mathbb{E}_1^3$. Finally, we find relationship between relatively normal-slant helices and slant helices on timelike as well as spacelike surfaces.

\section{Preliminaries}

First of all we give some brief introduction of Minkowski $3$-space $E_1^3$. The space  $\mathbb{E}_1^3$ is a three dimensional real vector space endowed with the dot product
\begin{equation}
\langle x, y \rangle = - x_1y_1 + x_2y_2 + x_3y_3,
\end{equation}
\n where $x=(x_1, x_2, x_3)$, $y=(y_1, y_2, y_3)$ $\in E_1^3$. This space is also known as Lorentz-Minkowski space. A vector $x \in E_1^3$ is said to be spacelike when $\langle x, x \rangle >$ 0  or x = 0, timelike when  $\langle x, x \rangle <$ 0 and lightlike(null) when $\langle x, x \rangle$ = 0. A curve  $\gamma: I \rightarrow E_1^3$ is called spacelike, timelike or lightlike when the velocity vector of the curve is spacelike, timelike or lightlike, respectively. While a surface $M$ is called spacelike, timelike or lightlike when the unit normal of the surface is timelike, spacelike or lightlike, respectively.

\n The Lorentzian cross product of $x$ = ($x_1, x_2, x_3$) and $y$ = ($y_1, y_2, y_3$) $\in E_1^3$ is defined as follows

$x \times y =$ $
\begin{vmatrix}
 e_1 & -e_2 & -e_3\\
 x_1 & x_2 & x_3\\
 y_1 & y_2 & y_3 
 \end{vmatrix}
$ = $(x_2y_3 - x_3y_2, x_1y_3 - x_3y_1, x_2y_1 - x_1y_2)$,

\n where $e_i = (\delta_{i1}, \delta_{i2}, \delta_{i3})$, $ \delta_{ij} $ is Kronecker delta and $e_1 \times e_2$ = $-e_3$, $e_2 \times e_3 = e_1$, $e_3 \times e_1 = -e_2$.

\n Let $\{ T, n, b\}$ be the moving Frenet frame along the curve $\gamma$ with arc-length parameter $s$ with curvature $\kappa$ and torsion $\tau$.

\n For a spacelike curve $\gamma$, the Frenet-Serret equations are given by

\begin{center}
$ \begin{bmatrix}
       
       T' \\
       n' \\
       b'
       \end{bmatrix}$ = $ \begin{bmatrix} 
                                      0 & \kappa & 0\\
                                      -\epsilon\kappa & 0 & \tau\\
                                      0 & \tau & 0
                                      \end{bmatrix}$ = $ \begin{bmatrix}
       
                                                                         T \\
                                                                          n\\
                                                                          b
                                                                  \end{bmatrix},$
                                                           
                                                           \end{center}       
\n where $\langle T, T \rangle$ = 1, $\langle n, n \rangle$ = $\epsilon$, $\langle b, b \rangle$ = $-\epsilon$, $\langle T, b \rangle$ = $\langle T, n \rangle$ = $\langle n, b \rangle$ = 0. When $\epsilon = 1$, $\gamma(s)$ is a spacelike curve with spacelike principal normal n and timelike binormal b while if $\epsilon $ = -1 then $\gamma$ is a spacelike curve with timelike principal normal n and spacelike binormal b.  
       
       \n For a timelike curve $\gamma$, the Frenet-Serret equations are given by
       
   \begin{center}    
       $ \begin{bmatrix}
       
        T' \\
       n' \\
       b'
       \end{bmatrix}$ = $ \begin{bmatrix} 
                                      0 & \kappa & 0\\
                                      \kappa & 0 & \tau\\
                                      0 & -\tau & 0
                                      \end{bmatrix}$ = $ \begin{bmatrix}
       
                                                                         T \\
                                                                          n\\
                                                                          b
                                                                  \end{bmatrix},$
       \end{center}
       
 \n where $\langle T, T\rangle$ = -1, $\langle n, n\rangle$ = 1, $\langle b, b\rangle$ = 1, $\langle T, n\rangle$ = $\langle T, b\rangle$ = $\langle n, b\rangle$ = 0.
       
\n \begin{definition}$[12]$ Let $v$ and $w$ be two spacelike vectors of $E_1^3$. Then we have the following:

\n (a) If $v$ and $w$ span a spacelike vector subspace of $E_1^3$, then there is a unique non-negative real number $\theta \geq 0$ such that $\langle v, w \rangle = \|v\|  \|w\| \cos\theta$.

\n (b) If $v$ and $w$ span a timelike vector subspace of $E_1^3$, then there is a unique non-negative real number $\theta \geq 0$ such that $\langle v, w \rangle = \|v\|  \|w\| \cosh\theta$.

\end{definition}

\begin{definition} $[12]$ Let $v$ be a spacelike vector and $w$ be a timelike vector in $\mathbb{E}_1^3$. Then there is a unique non-negative real number $\theta \geq$ 0, such that 
\begin{center}
\n $\langle v, w \rangle = \|v\|  \|w\| \sinh\theta$.

\end{center}
\end{definition}

\begin{definition} $[10]$ Let $v$ be a timelike vector and $w$ be a timelike vector in same time cone of $\mathbb{E}_1^3$, i.e. $\langle v, w\rangle < 0$. Then there is a unique non-negative real number $\theta \geq$ 0, such that 
\begin{center}
\n $\langle v, w \rangle = -\|v\|  \|w\| \cosh\theta$.

\end{center}
\end{definition}

\n Let $M$ be a smooth spacelike surface in $\mathbb{E}_1^3$ and $\gamma: I \rightarrow E_1^3$ be a unit speed spacelike curve on the surface. Then the Darboux frame $\{T, B = N \times T, N \}$ along the curve is well-defined and positively oriented along the curve, where $T$ is the tangent vector field of $\gamma$, $N$ is the unit normal of $M$ and $B$ is intrinsic normal of $\gamma$. The Darboux equations are given by
                                           
 \begin{align}\label{2.2}
 T^{'}   = \kappa_g B + \kappa_n N,\;
 B^{'}  = - \kappa_g T + \tau_g N,\;  
 N^{'}  =  \kappa_n T + \tau_g B, 
 \end{align}
 
\n where $\kappa_{g}$, $\kappa_{n}$ and $\tau_{g}$ are the geodesic curvature, normal curvature and geodesic torsion, respectively, and $\langle T, T\rangle = \langle B, B\rangle =1$, $\langle N, N\rangle = -1$, $\langle n, n\rangle = 1$.

\n By using \eqref{2.2} we get,

 \begin{equation}\label{2.3}
 \kappa^{2}  = \kappa_g^{2} -  \kappa_n^{2},\;
 \kappa_n  = \kappa \sinh \phi, \;
 \kappa_g  = \kappa \cosh \phi, \;
 \tau_g  = \tau + \phi^{'},
 \end{equation}
  \n where $\phi$ is the angle between the surface normal $N$ and the principal normal $n$ to the curve $\gamma$. 
  
  \n Now, let $M$ be a smooth timelike surface in $E_1^3$. Then the tangent space of the surface is timelike therefore a curve lying on the surface could be either timelike or spacelike. Let $\gamma : I \rightarrow E_1^3$ be a unit speed timelike curve lying on the surface $M$. Then Darboux equations are given by

\begin{equation}\label{2.4}
 T^{'}  = \kappa_g B + \kappa_n N, \;
 B^{'}  =  \kappa_g T - \tau_g N,  \;
 N^{'}  =  \kappa_n T + \tau_g B, 
 \end{equation}

 \n
where $\kappa_{g}$, $\kappa_{n}$ and $\tau_{g}$ are the geodesic curvature, normal curvature and geodesic torsion, respectively and $\langle T, T\rangle = -1, \langle B, B\rangle  = \langle N, N\rangle = \langle n, n\rangle = 1$.

 \n By using \eqref{2.4} we get,

 \begin{equation}\label{2.5}
 \kappa^{2}  = \kappa_g^{2} + \kappa_n^{2}, \;
 \kappa_n  = \kappa \cos \phi, \; 
 \kappa_g  = \kappa \sin \phi, \;
 \tau_g  = \tau + \phi^{'},
 \end{equation}
 
 \n where $\phi$ is the angle between the surface normal $N$ and the principal normal $n$ to the curve $\gamma$.
  
 \n
 If $\gamma : I \rightarrow E_1^3$ is a unit speed spacelike curve on the timelike surface $M$, then Darboux equations are given by
 
 \begin{equation}\label{2.6}
 T^{'}  = \kappa_g B - \kappa_n N, \;
 B^{'}  =  \kappa_g T + \tau_g N,  \;
 N^{'}  =  \kappa_n T + \tau_g B,
 \end{equation}
 
 \n
where $\kappa_{g}$, $\kappa_{n}$ and $\tau_{g}$ are the geodesic curvature, normal curvature and geodesic torsion, respectively, and $\langle T, T\rangle = \langle N, N\rangle  = \langle n, n\rangle = 1, \langle B, B\rangle = -1$.

 \n By using Eq.\eqref{2.6} we get,

 \begin{equation}\label{2.7}
 \kappa^{2}  = \kappa_n^{2} + \kappa_g^{2}, \;
 \kappa_n  = \kappa \cosh \phi, \; 
 \kappa_g  = \kappa \sinh \phi, \;
 \tau_g  = \tau + \phi^{'} ,
 \end{equation}
 
 \n where $\phi$ is the angle between the surface normal $N$ and the principal normal $n$ to the curve $\gamma$.

\begin{theorem} $[1]$ \label{2.4}
Let $\gamma$ be a unit speed spacelike curve in $E_1^3$. If the normal vector of the $\gamma$ is spacelike, then $\gamma$ is a slant helix if and only if one of the two functions 
$\frac{\kappa^2}{(\kappa^2 - \tau^2)^{3/2}} \left (\frac{\tau}{\kappa} \right )' $ \; and $ \frac{\kappa^2}{(\tau^2 - \kappa^2)^{3/2}} \left (\frac{\tau}{\kappa} \right )' $
is a constant everywhere $\tau^2 - \kappa^2$ does not vanish .
\end{theorem}
\begin{theorem} $[1]$ \label{2.5}
Let $\gamma$ be a unit speed timelike curve in $E_1^3$. Then $\gamma$ is a slant helix if and only if one of the two functions 
$\frac{\kappa^2}{(\kappa^2 - \tau^2)^{3/2}} \left (\frac{\tau}{\kappa} \right )' $ \; and $ \frac{\kappa^2}{(\tau^2 - \kappa^2)^{3/2}} \left (\frac{\tau}{\kappa} \right )' $
is a constant everywhere $\tau^2 - \kappa^2$ does not vanish .
\end{theorem}

\section{The axis of a spacelike relatively normal-slant helix on a spacelike surface}

\n Let $\gamma$ be a unit speed spacelike curve on an spacelike surface $M$ and $\{T, B, N\}$ be the Darboux frame along $\gamma(s)$. The curve $\gamma$ is called a relatively normal-slant helix if the vector field $B$ of $\gamma$ makes a constant angle with a fixed direction $d$. The unit vector $d$ is called the axis of the relatively normal-slant helix. In this section, we find the fixed vector(axis) of a spacelike relatively normal-slant helix via Darboux frame on a spacelike surface immersed in Minkowski $3$-space. We examine the two different cases of the axis $d$.

\n $\textbf{Case(1).}$
If the vector d is spacelike vector, then since the intrinsic normal B to the curve $\gamma$ is spacelike, from Definition 2.1(a), we have $\langle B, d\rangle = \cos \theta $ and from Definition 2.1(b), we have $\langle B, d\rangle = \cosh \beta $, where $\theta$ and $\beta$ are the fixed angles between the vectors $B$ and $d$, respectively.

\n $\textbf{(a)}$ Let $\langle B, d\rangle = \cosh \beta.$ If we differentiate this equation with respect to $s$ along the curve $\gamma$ and then by using \eqref{2.2}, we get $\langle B', d\rangle = 0\;\implies \langle -\kappa_g T + \tau_g N, d\rangle = 0 \;\implies \langle T, d\rangle = \frac{\tau_g}{\kappa_g} \langle N, d\rangle.$ Let us say $\langle N, d\rangle =$ a, then $\langle T, d\rangle = a\frac{\tau_g}{\kappa_g}$ and hence the axis can be written as $d = a\frac{\tau_g}{\kappa_g} T + \cosh \beta B - a N$.

\n Then
 $ \langle d, d\rangle$ = $a^2(\frac{\tau_g}{\kappa_g})^2 + \cosh^2 \beta - a^2 = 1$, which gives $a = \pm \frac{\kappa_g}{\sqrt{\kappa_g^2 - \tau_g^2}} \sinh \beta$.

\n Thus the spacelike axis $d$ is given as,
\begin{equation}\label{3.1}
d = \pm \frac{\tau_g}{\sqrt{\kappa_g^2 - \tau_g^2}} \sinh \beta T+ \cosh \beta B \mp \frac{\kappa_g}{\sqrt{\kappa_g^2 - \tau_g^2}} \sinh \beta N.
\end{equation}
If we differentiate $B'$ in \eqref{2.2} and $\langle B', d\rangle = 0$ with respect to $s$, we get 

\begin{equation} \label{3.2}
B''= (- \kappa_g' + \tau_g \kappa_n) T + (\tau_g^2 - \kappa_g^2) B + (\tau_g' - \kappa_g \kappa_n) N \; and \; \langle B'', d\rangle = 0.
\end{equation}
From \eqref{3.1} and \eqref{3.2}, we get
\begin{equation}
\langle B'', d\rangle = \pm \frac{(\kappa_g \tau_g' - \kappa_g' \tau_g)  - \kappa_n (\kappa_g^2 -\tau_g^2) }{\sqrt{\kappa_g^2 - \tau_g^2}} \sinh \beta -  (\kappa_g^2 -\tau_g^2) \cosh \beta = 0, \nonumber
\end{equation}
which implies

\begin{equation} \label{3.3}
\coth \beta = \pm  \left [ \frac{\kappa_g^2}{(\kappa_g^2 - \tau_g^2)^{3/2}} \left (\frac{\tau_g}{\kappa_g} \right )' - \frac{\kappa_n}{(\kappa_g^2 - \tau_g^2)^{1/2}} \right ]. 
\end{equation}

\n Now we prove that $d$ is a constant vector, i.e. $d'$ = 0. So we differentiate $d$ in (3.1), with respect to $s$ and obtain

\begin{align}
d'  =& \pm \sinh \beta \left [ (\frac{\tau_g}{\sqrt{\kappa_g^2 - \tau_g^2}}  )' T + \frac{\tau_g}{\sqrt{\kappa_g^2 - \tau_g^2}}(\kappa_g B + \kappa_n N)  \right ]  \nonumber \\
& \pm \sinh \beta \left [ -(\frac{\kappa_g}{\sqrt{\kappa_g^2 - \tau_g^2}} )' N - \frac{\kappa_g}{\sqrt{\kappa_g^2 - \tau_g^2}}(\kappa_n T + \tau_g B)  \right ] + \cosh \beta  ( -\kappa_g T + \tau_g N ), \nonumber
\end{align}
which implies
\begin{align}
d'  =& \pm \left (\sinh \beta \left [ (\frac{\tau_g}{\sqrt{\kappa_g^2 - \tau_g^2}}  )'  - \frac{\kappa_n \kappa_g}{\sqrt{\kappa_g^2 - \tau_g^2}} \right ] - \kappa_g \cosh \beta  \right ) T  \\
& \pm \left ( \sinh \beta \left [ \frac{\tau_g \kappa_n}{\sqrt{\kappa_g^2 - \tau_g^2}} - (\frac{\kappa_g}{\sqrt{\kappa_g^2 - \tau_g^2}} )' \right ] + \tau_g \cosh \beta \right ) N. \nonumber
\end{align}
Also from (3.3), we get
\begin{equation}
\cosh \beta = \pm \sinh \beta \; \frac{\tau_g' \kappa_g - \kappa_g' \tau_g - \kappa_n(\kappa_g^2 - \tau_g^2)}{(\kappa_g^2 - \tau_g^2)^{3/2}}. 
\end{equation}
Now from (3.4) and (3.5), we get

\begin{center}
$d' = \pm \sinh \beta$ $\begin{pmatrix}
                                   \frac{\tau_g'(\kappa_g^2 - \tau_g^2) - \tau_g \kappa_g \kappa_g' + \tau_g^2 \tau_g' - \kappa_n \kappa_g (\kappa_g^2 - \tau_g^2)}{(\kappa_g^2 - \tau_g^2)^{3/2}} \\
                                   - \frac{\tau_g' \kappa_g^2 - \kappa_g' \kappa_g \tau_g - \kappa_n \kappa_g(\kappa_g^2 - \tau_g^2)}{(\kappa_g^2 - \tau_g^2)^{3/2}}
                                  \end{pmatrix}$ T
                                  
                $\pm \sinh \beta$ $\begin{pmatrix}
                                   \frac{\tau_g \kappa_n(\kappa_g^2 - \tau_g^2) - \kappa_g' (\kappa_g^2 - \tau_g^2) + \kappa_g^2 \kappa_g' - \kappa_g \tau_g \tau_g'}{(\kappa_g^2 - \tau_g^2)^{3/2}} \\
                                   + \frac{\tau_g \tau_g' \kappa_g - \kappa_g' \tau_g^2 - \kappa_n \tau_g(\kappa_g^2 - \tau_g^2)}{(\kappa_g^2 - \tau_g^2)^{3/2}}
                                  \end{pmatrix}$ N.               
                                  
\end{center}
By straight forward calculation we see that the above expression vanishes, i.e. $d'$ = 0. Hence $d$ is a constant vector.

\n $\textbf{(b)}$  Let $\langle B, d\rangle = \cos \theta$. Proceeding in similar way as Case 1(a) we conclude that $\langle T, d \rangle = \frac{\tau_g}{\kappa_g} \langle N, d \rangle$. If we consider $\langle N, d\rangle = a$,  then $d = a \frac{\tau_g}{\kappa_g} T + \cos \theta B - a N$. Now since $d$ is spacelike, therefore 1 = $\langle d, d \rangle = a^2 (\frac{\tau_g}{\kappa_g})^2 + \cos^2\theta - a^2 $, so $ a = \pm \frac{\kappa_g}{\sqrt{\tau_g^2 - \kappa_g^2}} \sin \theta.$

\n
Thus the spacelike axis $d$ can be written as,
\begin{equation}\label{3.6}
d = \pm \frac{\tau_g}{\sqrt{\tau_g^2 - \kappa_g^2}} \sin \theta T + \cos \theta B \mp \frac{\kappa_g}{\sqrt{\tau_g^2 - \kappa_g^2}} \sin \theta N.
\end{equation}
Now, since $B'' = (- \kappa_g' + \tau_g \kappa_n) T + (\tau_g^2 - \kappa_g^2) B + (\tau_g' - \kappa_g \kappa_n) N$ and $\langle B'', d \rangle = 0$, proceeding as Case 1(a), we obtain
\begin{equation} \label{3.7}
\cot \theta = \mp  \left [ \frac{\kappa_g^2}{(\tau_g^2 - \kappa_g^2)^{3/2}} \left (\frac{\tau_g}{\kappa_g} \right )' - \frac{\kappa_n}{\tau_g^2 - \kappa_g^2)^{1/2}} \right ]. 
\end{equation}

\n
On the other hand we can do similar calculations to show that $d'$=0, consequently $d$ is a constant vector.

\n
$\textbf{Case 2.}$ Let the axis $d$ be timelike. Then by Definition 2.2, $\langle B, d\rangle = \sinh \alpha$, where $\alpha$ is the constant angle between $B$ and $d$. In this case, we get
\begin{equation}\label{3.8}
d = \pm \frac{\tau_g}{\sqrt{\kappa_g^2 - \tau_g^2}} \cosh \alpha \; T + \sinh \alpha\; B \mp \frac{\kappa_g}{\sqrt{\kappa_g^2 - \tau_g^2}} \cosh \alpha \;N,
\end{equation}
\begin{equation} \label{3.9}
\tanh \alpha = \pm  \left [ \frac{\kappa_g^2}{(\kappa_g^2 - \tau_g^2)^{3/2}} \left (\frac{\tau_g}{\kappa_g} \right )' - \frac{\kappa_n}{(\kappa_g^2 - \tau_g^2)^{1/2}} \right ]. 
\end{equation}

\section{The axis of a spacelike relatively normal-slant helix on a timelike surface}

Let $\gamma$ be a unit speed spacelike curve on a timelike surface $M$ and $\{T, B, N\}$ be the Darboux frame along $\gamma(s)$. In this section, we find the fixed vector (axis) of a spacelike relatively normal-slant helix via Darboux frame on a timelike surface immersed in Minkowski $3$-space. We examine the two different cases of the axis $d$. Since $T$ is spacelike, also $N$ is spacelike, B has to be timelike vector.

\n $\textbf{Case 1.}$
If $d$ is a timelike vector such that $B$ and $d$ are in same time-cone i.e. $\langle B, d\rangle < 0$. Then by the Definition 2.3, we have $\langle B, d\rangle = - \cosh \delta$, where  $\delta$ is the constant angle between $B$ and $d$. Differentiating $\langle B, d\rangle = - \cosh \delta$ with respect to $s$, we get $\langle B', d\rangle = 0$, then by using (2.6) for $B'$, we obtain $\langle T, d\rangle = - \frac{\tau_g}{\kappa_g} \langle N, d\rangle$. Let us say $  \langle N, d\rangle $= a, then the timelike axis $d$ can be written as
\begin{equation} \label{4.1}
d = - a \frac{\tau_g}{\kappa_g} T + \cosh \delta B + a N.
\end{equation}
Since $\langle d, d\rangle = -1$ therefore from (4.1) $a = \pm \frac{\kappa_g}{\sqrt{\kappa_g^2 + \tau_g^2}} \sinh \delta$ and hence

\begin{equation}\label{4.2}
d = \mp \frac{\tau_g}{\sqrt{\kappa_g^2 + \tau_g^2}} \sinh \delta\; T+ \cosh \delta\; B \pm \frac{\kappa_g}{\sqrt{\kappa_g^2 + \tau_g^2}} \sinh \delta \;N.
\end{equation}

\n Also $\langle B'', d \rangle = 0$, proceeding as Case 1(a) of the Section: 3 we obtain

\begin{equation} \label{4.3}
\coth \delta = \pm  \left [ \frac{\kappa_g^2}{(\kappa_g^2 + \tau_g^2)^{3/2}} \left (\frac{\tau_g}{\kappa_g} \right )' - \frac{\kappa_n}{(\kappa_g^2 + \tau_g^2)^{1/2}} \right ] .
\end{equation}

\n $\textbf{Case 2.}$
If $d$ is a spacelike vector, then by Definition 2.2, $\langle B, d\rangle =  \sinh \zeta$. In this case, we obtain

\begin{equation}\label{4.4}
d = \mp \frac{\tau_g}{\sqrt{\kappa_g^2 + \tau_g^2}} \cosh \zeta\; T - \sinh \zeta\; B \pm \frac{\kappa_g}{\sqrt{\kappa_g^2 + \tau_g^2}} \cosh \zeta \;N,
\end{equation}
\begin{equation} \label{4.5}
\tanh \zeta = \mp  \left [ \frac{\kappa_g^2}{(\kappa_g^2 + \tau_g^2)^{3/2}} \left (\frac{\tau_g}{\kappa_g} \right )' - \frac{\kappa_n}{(\kappa_g^2 + \tau_g^2)^{1/2}} \right ]. 
\end{equation}

\section{The axis of a timelike relatively normal-slant helix on a timelike surface}
\n  Let $\gamma$ be a unit speed timelike relatively normal slant helix on a timelike surface $M$. Then $T$ is a timellike vector, $N$ is a spacelike vector which imply that $B$ is a spacelike vector. In this situation, we have two different cases for the axis $d$ of the curve $\gamma$.

\n
$\textbf{Case 1.}$ If the vector d is spacelike vector, then from Definition: 2.1 (a), we have $\langle B, d\rangle = \cos \nu $ and from Definition: 2.1 (b), we have $\langle B, d\rangle = \cosh \xi $, where $\nu$ and $\xi$ are the constant angle between the vector $B$ and $d$ respectively.

\n \textbf{(a)} Let $\langle B, d\rangle = \cosh \xi$. Then differentiating the equation with respect to $s$ and using Darboux equation (2.4), we obtain
\begin{equation} \label{5.1}
d = - a \frac{\tau_g}{\kappa_g} T + \cosh \xi B + a N,
\end{equation}
where $a = \langle N, d\rangle$. Also $\langle d, d\rangle = 1$, which gives $a =  \pm \frac{\kappa_g}{\sqrt{\tau_g^2 - \kappa_g^2}} \sinh \xi$, and hence we obtain from (5.1)

\begin{equation}\label{5.2}
d = \mp \frac{\tau_g}{\sqrt{\tau_g^2 - \kappa_g^2}} \sinh \xi \; T + \cosh \xi \;B \pm \frac{\kappa_g}{\sqrt{\tau_g^2 - \kappa_g^2}} \sinh \xi \;N.
\end{equation}
Now by using $\langle B'', d\rangle = 0$, we obtain
\begin{equation} \label{5.3}
\coth \xi = \mp  \left [ \frac{\kappa_g^2}{(\tau_g^2 - \kappa_g^2)^{3/2}} \left (\frac{\tau_g}{\kappa_g} \right )' + \frac{\kappa_n}{(\tau_g^2 - \kappa_g^2)^{1/2}} \right ]. 
\end{equation}
\textbf{(b)} Let $\langle B, d\rangle = \cos \nu$. Then in similar way as Case 1(a), we obtain
 \begin{equation}\label{5.4}
d = \mp \frac{\tau_g}{\sqrt{\kappa_g^2 - \tau_g^2}} \sin \nu \; T + \cos \nu \; B \pm \frac{\kappa_g}{\sqrt{\kappa_g^2 + \tau_g^2}} \sin \nu \;N,
\end{equation}
\begin{equation} \label{5.5}
\cot \nu = \mp  \left [ \frac{\kappa_g^2}{(\kappa_g^2 - \tau_g^2)^{3/2}} \left (\frac{\tau_g}{\kappa_g} \right )' + \frac{\kappa_n}{(\kappa_g^2 - \tau_g^2)^{1/2}} \right ].
\end{equation}

\n $\textbf{Case 2.}$ When $d$ is timelike, by Definition 2.2 we have, $\langle B, d\rangle$ = $\sinh \psi$. In this case we obtain
\begin{equation}\label{5.6}
d = \mp \frac{\tau_g}{\sqrt{\tau_g^2 - \kappa_g^2}} \cosh \psi \; T + \sinh \psi \;B \pm \frac{\kappa_g}{\sqrt{\tau_g^2 - \kappa_g^2}} \cosh \psi \;N,
\end{equation}
\begin{equation} \label{5.7}
\tanh \psi = \mp  \left [ \frac{\kappa_g^2}{(\tau_g^2 - \kappa_g^2)^{3/2}} \left (\frac{\tau_g}{\kappa_g} \right )' + \frac{\kappa_n}{(\tau_g^2 - \kappa_g^2)^{1/2}} \right ]. 
\end{equation}

\section{Main Theorems}

In this section, we give main theorems that characterize relatively normal slant helices on a smooth surface immersed in Minkowski $3$-space.

\begin{theorem} \label{6.1}
A unit speed spacelike curve on a spacelike surface is a relatively normal slant helix if and only if any one of the following three functions, 

\n
(i) $\tanh \alpha = \mu(s) = \pm  \left [ \frac{\kappa_g^2}{(\kappa_g^2 - \tau_g^2)^{3/2}} \left (\frac{\tau_g}{\kappa_g} \right )' - \frac{\kappa_n}{(\kappa_g^2 - \tau_g^2)^{1/2}} \right ], $

\n
(ii) $ \coth \beta = \eta(s) = \pm  \left [ \frac{\kappa_g^2}{(\kappa_g^2 - \tau_g^2)^{3/2}} \left (\frac{\tau_g}{\kappa_g} \right )' - \frac{\kappa_n}{(\kappa_g^2 - \tau_g^2)^{1/2}} \right ], $

\n
(iii) $ \cot \theta = \psi(s) = \mp  \left [ \frac{\kappa_g^2}{(\tau_g^2 - \kappa_g^2)^{3/2}} \left (\frac{\tau_g}{\kappa_g} \right )' - \frac{\kappa_n}{(\tau_g^2 - \kappa_g^2)^{1/2}} \right ] $

\n
is a constant function.
\end{theorem}

\begin{proof}
Let the unit speed spacelike curve $\alpha$ on a surface $M$ with spacelike principal normal be a relatively normal slant helix, and hence the intrinsic normal $B$ to the curve makes a constant angle with a fixed direction. The relatively normal indicatrix i.e. the Gaussian mapping 
$B|_{\alpha} : I \rightarrow S_1^2$ along the curve $\gamma$ is a part of a circle on the Lorentzian unit sphere $S_1^2$. Therefore the relatively normal indicatrix has constant geodesic curvature and the normal curvature is one.

\n
For the curve $B|_{\alpha} : I \rightarrow S_1^2$, $B'$ is the tangent vector, and by using (2.2), we get

\begin{center}
$B' = -\kappa_g\; T + \tau_g\; N$,
\end{center}
\begin{center}
$B'' =  (- \kappa_g' + \tau_g \kappa_n) T + (\tau_g^2 - \kappa_g^2) B + (\tau_g' - \kappa_g \kappa_n) N.$
\end{center}
\n
Now, we have $N \times T = B, B\times N = T$ and $T\times B = N$. Thus

\begin{center}
$B' \times B'' = \tau_g(\kappa_g^2 - \tau_g^2) T + (\kappa_g \tau_g' -\tau_g  \kappa_g' + \kappa_n(\tau_g^2 - \kappa_g^2)) B - \kappa_g (\tau_g^2 - \kappa_g^2) N$,
\end{center}
which implies
\n
$\| B' \times B'' \|^2 = - (\kappa_g^2 - \tau_g^2)^3 + (\kappa_n (\tau_g^2 - \kappa_g^2) + \kappa_g^2 (\frac{\tau_g}{\kappa_g})')^2$ and $\| B' \|^2 = \kappa_g^2 - \tau_g^2$.

\n Let $\bar{\kappa}$ = Curvature of the relatively normal indicatrix. Then

\begin{center}
$\bar{\kappa} = \frac{\| B' \times B'' \|}{\| B' \|^3} = \sqrt{\sigma_B^2 - 1}$, 
\end{center}
where $\sigma_B = \pm \frac{1}{(\kappa_g^2 - \tau_g^2)^{\frac{3}{2}}} \left ( \kappa_g^2(\frac{\tau_g}{\kappa_g})' - \kappa_n(\kappa_g^2 - \tau_g^2) \right )$ and from (3.3), we have
\begin{equation} 
\coth \beta = \pm  \left [ \frac{\kappa_g^2}{(\kappa_g^2 - \tau_g^2)^{3/2}} \left (\frac{\tau_g}{\kappa_g} \right )' - \frac{\kappa_n}{(\kappa_g^2 - \tau_g^2)^{1/2}} \right ]. \nonumber
\end{equation}
Now if $\bar{\kappa_g}$ and $\bar{\kappa_n}$ are the geodesic and normal curvature of the relatively normal indicatrix, respectively. We have $\bar{\kappa}^2 = \bar{\kappa_g}^2 - \bar{\kappa_n}^2$ and we know that normal curvature of any spherical curve is 1, so $\bar{\kappa_n} = 1$, this imply that $\bar{\kappa_g}^2 - 1 = \sigma_B^2 - 1$, which implies
\begin{center}
$\bar{\kappa_g}^2 = \pm \sigma_B = \pm \coth \beta$.
\end{center}
Thus $\bar{\kappa_g}(s) $ is a constant function if and only if $\coth \beta(s)$ is a constant function. That means, relatively normal indicatrix of $\gamma$ is a part of a circle on $S_1^2$ if and only if 
\begin{equation} 
\coth \beta (s)= \pm  \left [ \frac{\kappa_g^2}{(\kappa_g^2 - \tau_g^2)^{3/2}} \left (\frac{\tau_g}{\kappa_g} \right )' - \frac{\kappa_n}{(\kappa_g^2 - \tau_g^2)^{1/2}} \right ] \nonumber
\end{equation}
is a constant function.

\n
Proof of (i) and (iii) can be done similarly.
\end{proof}

\begin{theorem} \label{6.2}
A unit speed spacelike curve on a timelike surface is a relatively normal slant helix if and only if any one of the following two functions, 

\n
(i) $\coth \delta = \lambda (s) = \pm  \left [ \frac{\kappa_g^2}{(\kappa_g^2 + \tau_g^2)^{3/2}} \left (\frac{\tau_g}{\kappa_g} \right )' - \frac{\kappa_n}{(\kappa_g^2 + \tau_g^2)^{1/2}} \right ]$, 

\n
(ii) $\tanh \zeta = \rho (s) = \mp  \left [ \frac{\kappa_g^2}{(\kappa_g^2 + \tau_g^2)^{3/2}} \left (\frac{\tau_g}{\kappa_g} \right )' - \frac{\kappa_n}{(\kappa_g^2 + \tau_g^2)^{1/2}} \right ]$

\n is a constant function.

\n
The proof is similar to the proof of Theorem 6.1 .
\end{theorem}

\begin{theorem} \label{6.3}
A unit speed timelike curve on a timelike surface is a relatively normal slant helix if and only if any one of the following three functions, 

\n (i) $\coth \xi = \chi (s) = \mp  \left [ \frac{\kappa_g^2}{(\tau_g^2 - \kappa_g^2)^{3/2}} \left (\frac{\tau_g}{\kappa_g} \right )' + \frac{\kappa_n}{(\tau_g^2 - \kappa_g^2)^{1/2}} \right ],$

\n (ii)  $\cot \nu = \omega (s) = \mp  \left [ \frac{\kappa_g^2}{(\kappa_g^2 - \tau_g^2)^{3/2}} \left (\frac{\tau_g}{\kappa_g} \right )' + \frac{\kappa_n}{(\kappa_g^2 - \tau_g^2)^{1/2}} \right ],$

\n (iii) $\tanh \psi = \mu (s) = \mp  \left [ \frac{\kappa_g^2}{(\tau_g^2 - \kappa_g^2)^{3/2}} \left (\frac{\tau_g}{\kappa_g} \right )' + \frac{\kappa_n}{(\tau_g^2 - \kappa_g^2)^{1/2}} \right ]$

\n is a constant function.

\n
The proof is similar to the proof of Theorem 6.1 .
\end{theorem}

\begin{proposition}
Let $\gamma$ be a unit speed spacelike relatively normal slant helix on a spacelike surface. Then

\n
(i) $\gamma$ is a asymptotic curve on the surface if and only if $\gamma$ is a slant helix with the spacelike axis
\begin{equation}
d = \pm \frac{\tau_g}{\sqrt{\kappa_g^2 - \tau_g^2}} \sinh \beta T+ \cosh \beta B \mp \frac{\kappa_g}{\sqrt{\kappa_g^2 - \tau_g^2}} \sinh \beta N, \nonumber
\end{equation}
(ii) $\gamma$ is a asymptotic curve on the surface if and only if $\gamma$ is a slant helix with the spacelike axis
\begin{equation}
d = \pm \frac{\tau_g}{\sqrt{\tau_g^2 - \kappa_g^2}} \sin \theta T + \cos \theta B \mp \frac{\kappa_g}{\sqrt{\tau_g^2 - \kappa_g^2}} \sin \theta N, \nonumber
\end{equation}
(iii) $\gamma$ is a asymptotic curve on the surface if and only if $\gamma$ is a slant helix with the timelike axis
\begin{equation}
 d = \pm \frac{\tau_g}{\sqrt{\kappa_g^2 - \tau_g^2}} \cosh \alpha \; T + \sinh \alpha\; B \mp \frac{\kappa_g}{\sqrt{\kappa_g^2 - \tau_g^2}} \cosh \alpha \;N,\nonumber
\end{equation}
(for the Case 1(a), 1(b) and Case 2, respectively).
\end{proposition}
\begin{proof}
(i) Since $\gamma$ is asymptotic, therefore $\kappa_n = 0 $. From (2.3) it follows that $\kappa_g = \kappa$ and $\kappa_n = \kappa \sinh \phi = 0$ which imply that $\phi = 0 \implies \phi' = 0$ so $\tau_g = \tau + \phi' \implies \tau_g = \tau$.

\n
Now by substituting $\kappa_g = \kappa, \kappa_n = 0$ and $\tau_g = \tau$ in (3.3), we obtain
\begin{equation} 
\eta (s) = \pm  \left [ \frac{\kappa^2}{(\kappa^2 - \tau^2)^{3/2}} \left (\frac{\tau}{\kappa} \right )' \right ] \nonumber
\end{equation}
is a constant function. Then by Theorem 2.4, $\gamma$ is a slant helix. From (3.1) we get the spacelike axis of the slant helix as 
\begin{equation}
d = \pm \frac{\tau}{\sqrt{\kappa^2 - \tau^2}} \sinh \beta T+ \cosh \beta B \mp \frac{\kappa}{\sqrt{\kappa^2 - \tau^2}} \sinh \beta N. 
\end{equation}
Conversely, let $\gamma$ be a slant helix with the spacelike axis $d$ as given in the above equation. Then since $\gamma$ is a relatively normal slant helix also, the fixed axis is 
\begin{equation}
d = \pm \frac{\tau_g}{\sqrt{\kappa_g^2 - \tau_g^2}} \sinh \beta T+ \cosh \beta B \mp \frac{\kappa_g}{\sqrt{\kappa_g^2 - \tau_g^2}} \sinh \beta N. \nonumber
\end{equation}
 Since both $B$ and $n$ are in the same plane and both make constant angle with the same fixed vector $d$, angle between $B$ and $n$ is also constant, i.e. $\phi = constant$ and hence $\phi'$ = 0. From (2.3), we get $\tau_g = \tau$ and $\kappa_g = \kappa \implies \kappa_n = 0$, that means $\gamma$ is a asymptotic curve.  
\n 
Also, (ii) and (iii) can be proved similarly.
\end{proof}
\begin{proposition}
Let $\gamma$ be a unit speed timelike relatively normal slant helix on a timelike surface. Then

\n
(i) $\gamma$ is a asymptotic curve on the surface if and only if $\gamma$ is a slant helix with the spacelike axis
\begin{equation}
d = \mp \frac{\tau_g}{\sqrt{\tau_g^2 - \kappa_g^2}} \sinh \xi \; T + \cosh \xi \;B \pm \frac{\kappa_g}{\sqrt{\tau_g^2 - \kappa_g^2}} \sinh \xi \;N, \nonumber
\end{equation}
(ii) $\gamma$ is a asymptotic curve on the surface if and only if $\gamma$ is a slant helix with the spacelike axis
\begin{equation}
d = \mp \frac{\tau_g}{\sqrt{\kappa_g^2 - \tau_g^2}} \sin \nu \; T + \cos \nu \; B \pm \frac{\kappa_g}{\sqrt{\kappa_g^2 + \tau_g^2}} \sin \nu \;N, \nonumber
\end{equation}
(iii) $\gamma$ is a asymptotic curve on the surface if and only if $\gamma$ is a slant helix with the timelike axis
\begin{equation}
d = \mp \frac{\tau_g}{\sqrt{\tau_g^2 - \kappa_g^2}} \cosh \psi \; T + \sinh \psi \;B \pm \frac{\kappa_g}{\sqrt{\tau_g^2 - \kappa_g^2}} \cosh \psi \;N,
\nonumber
\end{equation}
(for the Case 1(a), 1(b) and Case 2, respectively)
\end{proposition}
\begin{proof}
Proof is same as the proof of the Proposition 6.1.
\end{proof}
\begin{proposition}
If $\gamma$ is a unit speed (spacelike and timelike) relatively normal slant helix with(spacelike and timelike, respectively) axis $d$ on a (spacelike/timelike and timelike, respectively) surface then $\gamma$ can not be a line of curvature.
\end{proposition}
\begin{proof}
Suppose $\gamma$ is a line of curvature, then $\tau_g = 0$ and by Theorem 6.1, we have

\n
$ \coth \beta = \pm  \left [ \frac{\kappa_g^2}{(\kappa_g^2 - \tau_g^2)^{3/2}} \left (\frac{\tau_g}{\kappa_g} \right )' - \frac{\kappa_n}{(\kappa_g^2 - \tau_g^2)^{1/2}} \right ]. $

\n
If we put $\tau_g = 0, \coth \beta = \mp \frac{\kappa_n}{\kappa_g} = \mp \frac{\kappa \sinh \phi}{\kappa \cosh \phi} = \mp \tanh \phi$, which has no solution.

\n
So, $\gamma$ can not be a line of curvature.
Other cases can be proved similarly.
\end{proof}

\begin{proposition}
Let $\gamma$ be a unit speed spacelike relatively normal slant helix with timelike axis $d$. Then $\gamma$ is a plane curve provided $\gamma$ is a line of curvature on $M$.
\end{proposition}

\begin{proof} By Theorem 6.1, we have
\begin{equation} 
\tanh \alpha = \pm  \left [ \frac{\kappa_g^2}{(\kappa_g^2 - \tau_g^2)^{3/2}} \left (\frac{\tau_g}{\kappa_g} \right )' - \frac{\kappa_n}{(\kappa_g^2 - \tau_g^2)^{1/2}} \right ]  \nonumber
\end{equation}
is a constant function. Assume $\gamma$ is a line of curvature then $\tau_g = 0$, putting this into above equation we get $\tanh \alpha = \mp \frac{\kappa_n}{\kappa_g} = \mp \tanh \phi$. This implies that $\alpha = \pm \phi $ , and since $\alpha$ is constant therefore $\phi'$ = 0. Thus we get $\tau_g = \tau + \phi' = \tau$ and since $\tau_g = 0$, $\tau = 0$ as a result $\gamma$ is a planer curve.
\end{proof}

\begin{corollary}
 Let $\gamma$ be a spacelike relatively normal slant helix with spacelike axis on a timelike surface $M$, then $d$ cannot be orthogonal to the tangent line of $\gamma$.
 \end{corollary}
 \begin{proof} Suppose $\gamma$ be a unit speed spacelike relatively normal slant helix with spacelike axis  on the timelike surface $M$. Then by using (4.4), we have
 \begin{equation}
\langle T, d\rangle = \mp \frac{\tau_g}{\sqrt{\kappa_g^2 + \tau_g^2}} \cosh \zeta. \nonumber
 \end{equation}
 From Proposition 6.3, we have $\tau_g \neq 0$ also $\cosh \zeta $ never zero. This implies that $\langle T, d\rangle \neq 0$. Hence $T$ is not orthogonal to $d$.
  \end{proof}
 \begin{corollary}
  Let $\gamma$ be a spacelike relatively normal slant helix with timelike axis on a spacelike surface $M$, then $d$ cannot be orthogonal to $N$.
  \end{corollary}
  \begin{proof}
 Let $\gamma$ be a unit speed spacelike relatively normal slant helix with timelike axis  on the spacelike surface $M$. Then from (3.8), we have
\begin{equation}
 \langle N, d\rangle = \mp \frac{\kappa_g}{\sqrt{\kappa_g^2 - \tau_g^2}} \cosh \alpha. \nonumber
\end{equation}
 Since $\cosh \alpha $ never zero and $\kappa_g \neq 0 $ (unless $ \kappa = 0$) because $\kappa_g = \kappa \cosh \phi$ therefore $\langle N, d\rangle \neq 0$.
\end{proof}
\begin{corollary}
Let $\gamma$ be a timelike unit speed relatively normal slant helix with timelike axis d on a timelike surface $M$. Then

\n
(i) $d$ cannot be orthogonal to the tangent line of $\gamma$,

\n
(ii) $d$ is orthogonal to the surface normal $N$ if and only if $\gamma$ is a geodesic.
\end{corollary}
\begin{proof}
Let $\gamma$ be a timelike unit speed relatively normal slant helix with timelike axis d on a timelike surface $M$. Then by (5.6), we have

\n
 
\begin{equation}
\langle T, d\rangle = \mp \frac{\tau_g}{\sqrt{\tau_g^2 - \kappa_g^2}} \cosh \psi.\nonumber
\end{equation}
From Proposition 6.3, we have $\tau_g \neq 0$ and we know $\cosh \psi$ never zero, this imply that $\langle T, d\rangle \neq 0$, so $d$ cannot be orthogonal to $T$.

\n
Also from (5.6), we have 
\begin{equation}
\langle N, d\rangle = \pm \frac{\kappa_g}{\sqrt{\tau_g^2 - \kappa_g^2}} \cosh \psi. \nonumber
\end{equation}
Thus $ \langle N, d\rangle = 0$ iff $\kappa_g$ = 0. Hence $d$ is orthogonal to $N$ if and only if $\gamma$ is a geodesic curve.
\end{proof}

\noindent\author{Akhilesh Yadav}\\
\date{Department of Mathematics, Institute of Science, \\Banaras Hindu University, Varanasi-221005, India}\\
\maketitle {\noindent E-mail: akhilesha68@gmail.com}

\noindent\author{Ajay Kumar Yadav}\\
\date{Department of Mathematics, Institute of Science, \\Banaras Hindu University, Varanasi-221005, India}\\
\maketitle {\noindent E-mail: ajaykumar74088@gmail.com}

\end{document}